\documentclass[11pt,reqno]{amsart}

\usepackage{graphicx}
\usepackage{indentfirst,csquotes}

\usepackage{amssymb,amsthm,amsmath}
\usepackage{xcolor,paralist,hyperref,titlesec,fancyhdr,etoolbox}

\titleformat{\section}[display]{\normalfont\huge\bfseries\centering}{\centering\chaptertitlename\thechapter}{10pt}{\Large}
\titlespacing*{\section}{0pt}{0ex}{0ex}

\hypersetup{ colorlinks=true, linkcolor=black, filecolor=black, urlcolor=black }

\usepackage{lipsum}
\usepackage{booktabs}

\begin{document}
\title{Integration by differentiation} 
\author[]{Alexandre M. Afonso}
\maketitle

\let\thefootnote\relax
\footnotetext{MSC2020: Primary 00A05, Secondary 34-02, 34A05, 97I40, 26A33} 

\begin{abstract}
In this note, we address a traditional pedagogical asymmetry between algorithmic differentiation and heuristic integration by introducing an extra parametric-order differential framework. By treating the order of differentiation as a parametric axis  (discrete or continuous), classical integration is reframed as a local differential evaluation at the non-positive boundary level. This analytical shift alleviates some traditional routines of integration by parts, by replacing the quadrature of a product of elementary functions with a direct differential series that either truncates into a finite algebraic sum or resolves through cyclic symmetries. This parametric-order approach also bridges the gap with fractional calculus, showing that non-integer derivatives emerge as a seamless extension of the same operational dimension. By reducing global quadrature to local differentiation, this note offers an efficient, conceptual tool for undergraduate analysis.
\end{abstract} 


\bigskip

\section*{Introduction}

Every student of calculus encounters a striking asymmetry: differentiation is a straightforward, local operation governed by deterministic rules, whereas integration remains almost an art of global heuristic patterns of inverse operations. We teach one mechanically, yet, when reversing the process, we lead students to abandon their algorithmic safety and rely on intuition, clever substitutions, or memorized tabulations.

This traditional approach treats the derivative derivative and the integral as separate operations. In this note, we explore this historical dichotomy, by viewing the order of the mathematical operation as a parameter (discrete or continuous), we map the behavior of elementary functions along a single, parametric-order axis. 

The analytical and pedagogical rewards are immediate. The symmetric notation allows that a single parameter to govern both rates of change and accumulations, causing the routines of multiple integration by parts to be replaced by a direct differential series, harmonic structures resolve through phase symmetries and fractional derivatives emerge naturally as coordinate evaluations. The global quadrature collapses into a local procedure, simplifying undergraduate analysis.

\section*{The Core Idea and Notation}
Let $f(x)$ be a well-defined elementary function whose integral and differential operations are analytically mapped through a parametric-order  map. This operator is indexed by a discrete parameter $n \in \mathbb{N}$ or, more generally, a  continuous variable $s \in \mathbb{R}$. We define the generalized state of a function $f(x)$ at order level $n$ using the notation $f(x, n)$, bounded by 
the following foundational axioms:
\begin{equation}
\begin{aligned}
    f(x, 1) &\equiv f(x) && \text{(Level 1 represents the function itself)} \\
    f(x, 0) &\equiv \int f(x) \, dx && \text{(Level 0 represents the primitive)}
    \label{1}
\end{aligned}
\end{equation}

By linear induction originating from the Fundamental Theorem of Calculus~\cite{rudin1976}, shifting the index order $n$ dynamically maps the successive differential and integral operators. The displacement of the parameter follows these symmetric translation rules:
\begin{equation}
\begin{aligned}
     \frac{d}{d x} \left[ f(x, n) \right] &= f(x, n+1) \quad \text{(Differentiation elevates the order)} \\
    \int f(x, n) \, dx &= f(x, n-1) \quad \text{(Integration reduces the order)}
\end{aligned}
\end{equation}

The core idea in this work lies in obtaining a general analytical expression for $f(x, n)$, established from the baseline state $f(x, 1)$, and compute the classical primitive simply by evaluating the parameter-order dimension at the zero level ($n = 0$), using identity axiom of equation \eqref{1}: $f(x, 0) \equiv \int f(x)$.

\subsection*{Example: $f(x,1) = x^3$.}
To understand how the generator $f(x,n)$ is constructed, we examine the inductive mechanism using a concrete example: $f(x,1) = x^3$. By tracking the successive standard derivatives of this cubic function, we systematically elevate the differential index level by level:
\begin{align*}
    f(x, 1) = x^3 , \quad
    f(x, 2) = 3x^2 , \quad
    f(x, 3) = 6x , \quad 
    f(x, 4) = 6  , \quad
    f(x, 5) = 0 \dots   
\end{align*}
with an arithmetic pattern emerging when we re-express the coefficients using factorials, and preserving the relationship between the starting exponent and the order level $n$,
\begin{align*}
    f(x,1) = \frac{3!}{3!}x^3,\; f(x,2) = \frac{3!}{2!}x^2, \; f(x,3) = \frac{3!}{1!}x^1, \; f(x,4) = \frac{3!}{0!}x^0, \; f(x,5) = 0 \dots
\end{align*}

By induction, this progression establishes the generator function for the cubic polynomial: 
\begin{equation}
f(x,n) = \frac{3!}{(4-n)!} x^{4-n}
\end{equation}
and the quadrature can be calculated by algebraic simplification at $n = 0$:
\begin{equation}
    \int x^3 \, dx \equiv f(x,0) = \frac{3!}{(3-0+1)!} x^{3-0+1} = \frac{3!}{4!} x^4 = \frac{3!}{4 \times 3!} x^4 = \frac{1}{4} x^4
\end{equation}
converging to the standard calculus result. Note that the vanishing terms within the higher-order derivatives in the sequence are natively absorbed by the negative factorial in the denominator. For any derivative state $n \geq 5$, the term $(4-n)!$ yields the factorial of a negative integer, which diverges to infinity and passively forces the entire reciprocal coefficient to zero.
 
To provide a complete overview of the general idea, Table~\ref{tab:universal_generators} summarizes the algebraic differential generators derived for some classical functions as defined at level $n=1$. 

\section*{The quadrature of function products}

We now extend this approach to the product of two fundamental functions, $f(x,1)$ and $g(x,1)$. The main goal is to determine the joint zero-level order, $fg(x,0)$, which maps directly to the classical indefinite integral:
\begin{equation}
    I = \int f(x, 1) \cdot g(x, 1) \, dx \equiv fg(x,0)
\end{equation}

By mapping the classical integration by parts framework~\cite{spivak2008} ($\int u \, dv = u v - \int v \, du$) onto the differential order dimension, we assign the baseline variables \(u = f(x,1)\) and \(dv = g(x,1) \, dx\). This yields the operational shifts \(du = f(x,2) \, dx\) and \(v = g(x,0)\). Successive iterations of this algebraic translation generate a structural chain of terms, shifting the differential indices in opposite directions. This dynamic recurrence establishes the following arithmetic ladder:
\begin{equation}
    I = f(x, 1)g(x, 0) - f(x, 2)g(x, -1) + f(x, 3)g(x, -2) - f(x, 4)g(x, -3) + \dots
\end{equation}

By tracing the pattern within this algebraic expansion, the entire infinite chain can be compactly mapped into a single differential series:
\begin{equation}
    \int f(x, 1) \cdot g(x, 1) \, dx = \sum_{j=0}^{\infty} (-1)^j \cdot f(x, 1+j) \cdot g(x, -j)
    \label{prod}
\end{equation}
where the signs alternate via $(-1)^j$, the order level of $f$ undergoes unit increments via $(1+j)$, and the order level of $g$ undergoes unit decrements via $(-j)$. 

This fundamental identity demonstrates that product quadrature can be reduced to direct algebraic substitution, provided the analytical generators $f(x, n)$ and $g(x, n)$ are established, such as the examples presented in Table~\ref{tab:universal_generators}. In practice, this infinite series is analytically applied at a finite step $k$ whenever the remaining residual vanishes to zero, either due to active polynomial termination or due to a cyclic loop closure. The convergence of this differential series is identical to the domain of success of classical integration by parts, operating effectively where traditional closed forms exist.

\subsection*{Example: $f(x,1)g(x,1) = x^2e^x$.}

To evaluate the quadrature of $f(x,1)g(x,1) = x^2e^x$, we assign the initial functions as $f(x,1) = x^2$ and $g(x,1) = e^x$. Using their respective polynomial and exponential generators from Table~\ref{tab:universal_generators} into the product series equation~\eqref{prod}, we obtain:
\begin{equation}
\begin{aligned}
    \int x^2 e^x \, dx &= \sum_{j=0}^{\infty} (-1)^j \cdot \left[ \frac{2!}{(2-j)!} x^{2-j} \right] \cdot e^x = 2! \cdot e^x \sum_{j=0}^{\infty} \frac{(-1)^j}{(2-j)!} \cdot x^{2-j} \\
     &= 2 \cdot e^x \left[ \frac{(-1)^0}{2!} x^2 + \frac{(-1)^1}{1!} x^1 + \frac{(-1)^2}{0!} x^0 \right] = e^x \left( x^2 - 2x + 2 \right)
\end{aligned}    
\end{equation}
which recovers the classic analytical primitive, bypassing traditional integration rules. Note that, again, an algebraic truncation occurs for all terms where $j \geq 3$, a region where the discrete denominator $(2-j)!$ hits the factorial of negative integers. 

\section*{The Fractional Derivative Space}

Transitioning from the discrete order parameter $n$ to a continuous variable $s \in \mathbb{R}$ also preserves the structural symmetry. Within the continuous differential space, the parameter $s$ anchors in the same fundamental states: $s=1$ denotes the native initial function, $s>1$ maps the continuous differential domain, and the origin $s=0$ extracts the classical primitive via direct evaluation, also verifying the identity axiom of equation~\eqref{1}, $f(x,0) \equiv \int f(x) \, dx$. Every elementary function is thus generalized as a continuous derivative axis $f(x,s)$, extending local calculus to an analytical manifold where any fractional order is evaluated by direct coordinate substitution.

To generalize the discrete order level into the continuous spectrum, we use the classical idea within the domain of fractional calculus~\cite{podlubny1999}, and replace the discrete factorial operator, $n!$, with its definitive analytic continuation, the Euler Gamma function, defined for $\Re(z) > 0$ as:
\begin{equation}
    \Gamma(z) = \int_0^\infty t^{z-1}e^{-t} \, dt
\end{equation}
which satisfies the fundamental factorial relation $\Gamma(z+1) = z!$. By shifting the parameter along this continuous axis, the entire spectrum of elementary functions can be generalized into the fractional state domain. Table~\ref{tab:universal_generators} also summarizes the fractional generators for these classical functions.

\begin{table}[htbp]
\centering
\caption{Function generators for discrete $f(x,n)$ and fractional $f(x,s)$ states.}
\label{tab:universal_generators}
\begin{tabular}{llll}
\toprule
$f(x,1)$ & \textbf{Discrete Generator } $f(x,n)$ & \textbf{Fractional Generator } $f(x,s)$ \\
\midrule
$x^p$ & $\dfrac{p!}{(p-n+1)!} x^{p-n+1} \quad (p, n \in \mathbb{N})$  & $\dfrac{\Gamma(p+1)}{\Gamma(p+2-s)} x^{p+1-s} \quad (p, s \in \mathbb{R})$ \\
\addlinespace
$e^{\lambda x}$ & $\lambda^{n-1} e^{\lambda x}$ & $\lambda^{s-1} e^{\lambda x}$ \\
\addlinespace
$\sin(x)$ & $\sin\left(x + (n-1)\dfrac{\pi}{2}\right)$ & $\sin\left(x + (s-1)\dfrac{\pi}{2}\right)$ \\
\addlinespace
$\cos(x)$ & $\cos\left(x + (n-1)\dfrac{\pi}{2}\right)$ & $\cos\left(x + (s-1)\dfrac{\pi}{2}\right)$ \\
\bottomrule
\end{tabular}
\end{table}

\subsection*{Example: The Semi-Integral of $f(x,1) = x^{\pi}$ at $s = 1/2$}

To demonstrate the precision of this approach, we compute the semi-integral of the irrational power function $f(x,1) = x^{\pi}$ and compare it with the direct evaluation using the classical fractional operator $D^{-1/2}$~\cite{podlubny1999}.

\subsubsection*{Resolution via the fractional algebric generator}
By evaluating the power function fractional generator in Table~\ref{tab:universal_generators} directly at the fractional coordinate $s = 1/2$ with a spatial exponent $p = \pi$, we obtain:
\begin{equation}
    f\left(x, s=\frac{1}{2}\right) = \frac{\Gamma(\pi+1)}{\Gamma\left(\pi + 2 - \frac{1}{2}\right)} x^{\pi + 1 - \frac{1}{2}}  \approx \frac{7.1878}{14.6797} x^{\pi + 0.5} \approx 0.4896 \cdot x^{\pi + 0.5}
\end{equation}
with the arguments $(\pi+1)$ and $(\pi+1.5)$ for the Euler Gamma function being real and strictly positive.

\subsubsection*{Comparison with the Classical Riemann-Liouville Convolution}
In traditional fractional calculus, computing the semi-integral of order $\alpha = 1/2$ requires fixing a geometric lower bound at the origin ($a=0$)~\cite{podlubny1999}. The Riemann-Liouville convolution operator is parameterized as follows~\cite{podlubny1999}:
\begin{equation}
    _{0}I^{1/2}_x [x^{\pi}] = \frac{1}{\Gamma(1/2)} \int_{0}^{x} \frac{t^{\pi}}{(x-t)^{1/2}} \, dt
\end{equation}

Resolving this singular integral on the real line requires a change of variables to invoke the classical Euler Beta function, $B(x,y) = \Gamma(x)\Gamma(y)/\Gamma(x+y)$~\cite{podlubny1999}, and obtain the global quadrature into a product of Gamma functions:
\begin{equation}
\begin{aligned}
    _{0}I^{1/2}_x [x^{\pi}] &= \frac{1}{\sqrt{\pi}} \cdot x^{\pi + 0.5} \cdot B\left(\pi+1, \frac{1}{2}\right) = \frac{1}{\sqrt{\pi}} \cdot x^{\pi + 0.5} \cdot \frac{\Gamma(\pi+1)\Gamma(1/2)}{\Gamma\left(\pi + 1 + \frac{1}{2}\right)} \\
    &= \frac{\Gamma(\pi+1)}{\Gamma\left(\pi + \frac{3}{2}\right)} x^{\pi + 0.5} \approx 0.4896 \cdot x^{\pi + 0.5}    
\end{aligned}
\end{equation}
since  $\Gamma(1/2) = \sqrt{\pi}$. This solution collapses into the same algebraic representation obtained via the continuous generator. The present methodology bypasses the need for global integration convolution routines, resolving fractional systems via simple local coordinate evaluations.

\section*{Discussion and Concluding Remarks}

The work presented in this note shows that by using a differential order to operate along a parametric integer/continuous axis, global quadrature can be reduced strictly to a localized coordinate evaluation. Properties such as polynomial truncation and cyclic harmonic loop closures emerge naturally from this single analytical domain. 

The conditions of applicability of this parametric order space are the same conditions of classical integration methodologies and the solutions obtained by the direct differential series are effective precisely within the domains where traditional integration by parts or substitutions resolve into stable closed forms. The existence of a universal operator generator function requires that the original functions possess stable analytic continuations into negative differential orders. 

Beyond its analytical utility, this parametric-order approach offers an immediate pedagogical advantage for undergraduate analysis. It relates standard primitives, integration by parts, and fractional calculus under a single parametric operator space. This allows educators to present these aspects of calculus not as isolated techniques, but as a continuous and interconnected spectrum.

\end{document}